\documentclass[9pt, a4paper]{amsart}
\usepackage{geometry} 
\usepackage{graphicx}
\usepackage[colorlinks, citecolor=blue, linkcolor=black]{hyperref}

\title[Two flat structures on minimal surfaces]{Two flat structures on minimal surfaces}
\author{Hojoo Lee}
\email{momentmaplee@gmail.com}

\theoremstyle{definition}

\numberwithin{equation}{section}

 \usepackage{float}

\begin{document}

\maketitle
 
In this expository article, we illustrate how two  \emph{flat} structures on minimal surfaces induce a harmonic function, which captures 
 the uniqueness of Enneper's algebraic surface with total Gauss curvature $-4\pi$. 
 
   \begin{enumerate}
%  \item[]  
  \item   Let  $\Sigma$ be a minimal surface in ${\mathbb{R}}^{3}$ with Gauss curvature $\mathcal{K}={\mathcal{K}}_{ {\mathbf{g}}_{{}_{\Sigma}}} \leq 0$. 
  Given the unit normal vector field $\mathbf{N}$ on  $\Sigma$ and a constant unit vector field ${\mathbf{V}}(p)={\mathbf{V}}$ in ${\mathbb{R}}^{3}$, introduce the 
  angle function 
\begin{center}
$
\mathbf{N}_{\mathbf{V}} (p):= {\langle \mathbf{N}(p), {\mathbf{V}} \rangle}_{{\mathbb{R}}^{3}}, \quad p \in \Sigma.
$
\end{center}
  \item \textbf{Ricci's flat structure.} The do Carmo--Peng proof  \cite{CP1979} of generalized \emph{Bernstein Theorem} that planes are the 
  only stable minimal surfaces in ${\mathbb{R}}^{3}$ uses the Ricci condition \cite[Equation (2.9)]{CP1979} 
\begin{center}
$
4  \mathcal{K} =  {\triangle}_{ {\mathbf{g}}_{{}_{\Sigma}} } \,  \ln \left( -   \mathcal{K}  \right),
$
\end{center}
which geometrically means Ricci's Theorem \cite{BM2017, Law1970, Law1971} that the conformally changed metric
\begin{center}
$
 {\mathcal{G}}_{{\text{Ricci}}} = {\left( -  {\mathcal{K}}  \right)}^{\frac{1}{2}}  {\mathbf{g}_{{}_{\Sigma}}} 
$
\end{center}
is \emph{flat} when $\mathcal{K}<0$. Ricci \cite[p. 124]{Bl1950} discovered that that every metric satisfying Ricci condition can be realized on a minimal 
surfaces in ${\mathbb{R}}^{3}$. See Lawson's counterexample \cite[Remark 12.1]{Law1970} for metics having a \emph{flat} point, and also his higher codimensional generalization 
 \cite{Law1971} of Ricci Theorem. 
    \item \textbf{Chern's flat structure.} Chern's non-complex-analytic proof \cite{Chern1969} of  \emph{Bernstein Theorem} that planes are the only entire minimal graphs in ${\mathbb{R}}^{3}$ uses the intriguing equality (with $\mathbf{V}=(0, 0, 1)$): 
\begin{center}
$
 \mathcal{K} =  {\triangle}_{ {\mathbf{g}}_{{}_{\Sigma}} } \,  \ln \left( 1+ \mathbf{N}_{\mathbf{V}}  \right), 
$
\end{center}
which geometrically means \cite[Equation (3)]{Chern1969} that the following conformally changed metric
\begin{center}
$
 {\mathcal{G}}_{{\text{Chern}}} = {\left(  1+ \mathbf{N}_{\mathbf{V}}  \right)}^{2}   {\mathbf{g}}_{{}_{\Sigma}}
$
\end{center}
is \emph{flat} when $\mathbf{N}_{\mathbf{V}}>-1$. (In \cite{Chern1969}, $\mathbf{N}_{(0, 0, 1)} = \frac{1}{\,\sqrt{\,1+ {f_{x}}^{2} + {f_{y}}^{2} \,}\,}>0$ on the minimal graph $z=f(x, y)$.)
 \item \textbf{Induced harmonic function.} When ${\mathcal{K}}<0$ and $\mathbf{N}_{\mathbf{V}} > -1$, since metrics $ \widetilde{g} ={\mathcal{G}}_{\text{Ricci}} $ and 
${\mathcal{G}}_{\text{Chern}} = \star \widetilde{g}  =  \frac{  { \left( 1 + {\mathcal{N}}_{\mathbf{V}}  \right)}^{2}  }{  {\left( - \mathcal{K}  \right)}^{\frac{1}{2}} }  {\mathcal{G}}_{\text{Ricci}}$ are \emph{flat}, by the curvature formula of conformally changed metric
\begin{center}
$
0=  {\mathcal{K}}_{ \star \widetilde{g} } = \frac{1}{\star}   {\mathcal{K}}_{  \widetilde{g} }  - \frac{1}{2 \star} {\triangle}_{ \widetilde{g} } \ln \star 
=0  - \frac{1}{2 \star} {\triangle}_{ \widetilde{g} } \ln \star,
$
\end{center}
we find the harmonicity (with respect to metrics ${\mathcal{G}}_{\text{Ricci}}$, ${\mathcal{G}}_{\text{Chern}}$, $ {\mathbf{g}_{{}_{\Sigma}}} $) of our Chern-Ricci function
\begin{center}
$
   \ln  \star  =   \ln  \left( \frac{  { \left( 1 + {\mathcal{N}}_{\mathbf{V}}  \right)}^{2}  }{  {\left( - \mathcal{K}  \right)}^{\frac{1}{2}} } \right).
$
\end{center}
 \item \textbf{Chern, Ricci, and Enneper \cite[Theorem 3.1]{Lee2017}.} If a minimal surface in ${\mathbb{R}}^{3}$ has a constant Chern-Ricci function, then it 
 is (a part of) an Enneper surface, up to  isometries and homotheties. 
  \end{enumerate}

 Enneper's surfaces have no \emph{flat} points. We list more uniqueness results. Bernstein-Mettler Theorem \cite{BM2017} characterizes the Enneper surface as 
 the \textbf{nonflat} minimal surface with vanishing entropy differential. P\'{e}rez \cite{Perez2007} gave a novel variational characterization of half of the Enneper surface containing 
 a straight line. White's original remarkable conjecture \cite{White 1996} states that \emph{the half Enneper surface is the unique properly embedded \textbf{nonflat} orientable 
 area minimizing surface with line boundary and quadratic area growth.}

\begin{figure}[H]
 \centering
 \includegraphics[height=4.50cm]{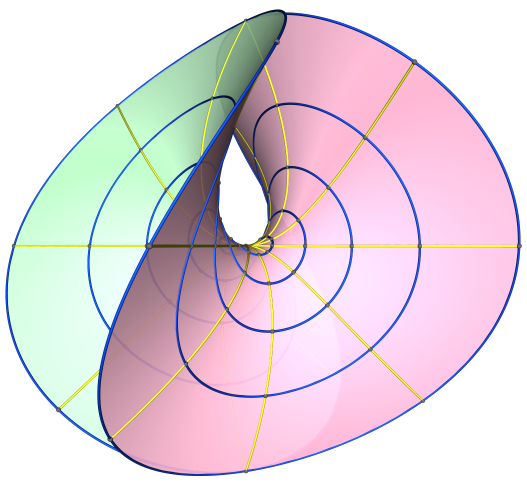}
 \caption{\small{An approximation \cite{WebEnn} of a part of Enneper's surface with total curvature $-4 \pi$}} 
 \end{figure}

 \bigskip

\end{document}